\documentclass[12pt,reqno,a4wide]{amsart}

\usepackage{tikz} 
\usepackage{calrsfs}
\usepackage{mathrsfs}

\usepackage{array}
\usepackage{hyperref}

\usetikzlibrary{decorations.pathreplacing}
\usetikzlibrary{fit}

\usepackage{pgfplots}

\allowdisplaybreaks

\catcode`,\active

\catcode`\,12

\begin{document}
\bibliographystyle{plain}

%
%

	\title
	{The amplitude of Motzkin paths}

	\author[H. Prodinger ]{Helmut Prodinger }
	\address{Department of Mathematics, University of Stellenbosch 7602, Stellenbosch, South Africa
	and
	Department of Mathematics and Mathematical Statistics,
	Umea University,
	907 36 Umea, 	 Sweden  }
	\email{hproding@sun.ac.za}

	
	\begin{abstract}
The amplitude of Motzkin paths was recently introduced, which is basically twice the height. We analyze this parameter
using generating functions.
	\end{abstract}
	

\maketitle

\section{Introduction}

A Motzkin path consists of up-steps, down-steps, and horizontal steps, see sequence A091965 in \cite{OEIS} and the references given there. As Dyck paths, they start at the origin and end, after $n$ steps 
again at the $x$-axis, but are not allowed to go below the $x$-axis. The height of a Motzkin path is the highest $y$-coordinate that the path
reaches. The average height of such random paths of length $n$ was considered in an early paper \cite{stack}, it is asymptotically given by
$\sqrt{\frac{\pi n}{3}}$. 

In the recent paper \cite{irene} an interesting new concept was introduced: the \emph{amplitude}. It is basically twice the height, but with a twist.
If there exists a horizontal  step on level $h$, which is the height, the amplitude is $2h+1$, otherwise it is $2h$. To clarify the concept, we  
created a list of all 9 Motzkin paths of length 4 with height and amplitude given.

\begin{center}
	
	\begin{table}[h]
		\begin{tabular}{c | c | c  |c}
			\text{Motzkin path  }&\text{horizontal   on maximal level}&\text{height}&\text{amplitude}\\
			 	\hline\hline
				\begin{tikzpicture}[scale=0.4]
							
						\draw[ultra thick] (0,0) to (1,0) to (2,0) to (3,0) to (4,0)  ;
									\end{tikzpicture}
				 & \text{Yes}& 0& 1\\
			
			\hline
			 
			\begin{tikzpicture}[scale=0.4]
				
				\draw[ ultra thick] (0,0) to (1,0) to (2,0) to (3,1) to (4,0)  ;
			\end{tikzpicture}
			& \text{No}& 1& 2\\
			
			\hline
			 \begin{tikzpicture}[scale=0.4]
			 	
			 	\draw[ ultra thick] (0,0) to (1,1) to (2,0) to (3,1) to (4,0)  ;
			 \end{tikzpicture}
			 & \text{No}& 1& 2\\
			 
			 \hline
			\begin{tikzpicture}[scale=0.4]
				
				\draw[ ultra thick] (0,0) to (1,0) to (2,1) to (3,1) to (4,0)  ;
			\end{tikzpicture}
			& \text{Yes}& 1& 3\\
			
			\hline
			\begin{tikzpicture}[scale=0.4]
				
				\draw[ ultra thick] (0,0) to (1,1) to (2,1) to (3,1) to (4,0)  ;
			\end{tikzpicture}
			& \text{Yes}& 1& 3\\
			\hline
			\begin{tikzpicture}[scale=0.4]
				
				\draw[ ultra thick] (0,0) to (1,1) to (2,1) to (3,0) to (4,0)  ;
			\end{tikzpicture}
			& \text{Yes}& 1& 3\\
			\hline
			\begin{tikzpicture}[scale=0.4]
				
				\draw[ ultra thick] (0,0) to (1,1) to (2,0) to (3,0) to (4,0)  ;
			\end{tikzpicture}
			& \text{No}& 1& 2\\
			\hline 
			\begin{tikzpicture}[scale=0.4]
				
				\draw[ ultra thick] (0,0) to (1,0) to (2,1) to (3,0) to (4,0)  ;
			\end{tikzpicture}
			& \text{No}& 1& 2\\
			
			\hline
              \begin{tikzpicture}[scale=0.4]
				
				\draw[ ultra thick] (0,0) to (1,1) to (2,2) to (3,1) to (4,0)  ;
			\end{tikzpicture}
			& \text{No}& 2& 4\\
			\hline
		\end{tabular}

	\end{table}	
\end{center}

The goal of this paper is to investigate this new parameter; in the next section, generating functions will be given, in the following section explicit enumerations,
involving trinomial coefficients  $\binom{n,3}{k}=[t^k](1+t+t^2)^n$ (notation following Comtet's book \cite{comtet}).
In the last section, the intuitive result that the average amplitude is about twice the average height, is confirmed, and then
it will be shown, that, asymptotically, there are about as many Motzkin paths with/without horizontal steps on the maximal level.

\section{Generating functions}

Let 
\begin{equation*}
M^{\le h}(z)=\sum_{n\ge0}[\text{number of Motzkin paths of length $n$ and height $\le h$}]z^n.
\end{equation*}
For the computation, let $f_i=f_i(z)$ be the generating function of Motzkin-like paths, bounded by height $h$, but ending at level $i$.
Distinguishing the last step, we get
\begin{equation*}
f_0=1+zf_0+zf_1,\quad f_i=zf_{i-1}+zf_i+zf_{i+1}\ \ \text{for}\ 0<i<h, \quad f_h=zf_{h-1}+zf_h.
\end{equation*}
This system is best written in matrix form:
\begin{equation*}
\begin{pmatrix}
	1-z&-z&0&\dots\\
	-z&1-z&-z&0&\dots\\
	0&-z&1-z&-z&0&\dots\\
	\ddots&\ddots&\ddots\\
	\phantom{0}&\phantom{0}&\phantom{0}&\phantom{0}&-z&1-z
\end{pmatrix}
\begin{pmatrix}
f_0\\f_1\\f_2\\ \vdots\\ f_h
\end{pmatrix}
=\begin{pmatrix}
	1\\0\\0\\ \vdots\\ 0
\end{pmatrix}
\end{equation*}
Let $D_n$ be the determinant of the system matrix, with $n$ rows and columns. Using Cramer's rule to solve a linear system, one finds
\begin{equation*}
M^{\le h}(z)=f_0=\frac{D_h}{D_{h+1}}.
\end{equation*}
Expanding the determinant along the first row, we get the recursion $D_n=(1-z)D_{n-1}-z^2D_{n-2}$, and $D_1=1-z$, $D_0=1$. Solving,
\begin{equation*}
D_n=\frac1{\sqrt{1-2z-3z^2}}\bigg[\biggl(\frac{1-z+\sqrt{1-2z-3z^2}}{2}\biggr)^{n+1}-\biggl(\frac{1-z-\sqrt{1-2z-3z^2}}{2}\biggr)^{n+1}\bigg].
\end{equation*}
If one deals with enumeration of Motzkin-like objects, the substitution $z=\frac{v}{1+v+v^2}$ makes the expressions prettier:
\begin{equation*}
D_n=\frac1{1-v^2}\frac{1-v^{2n+2}}{(1+v+v^2)^n}
\end{equation*}
and further
\begin{equation*}
M^{\le h}(z)=f_0=\frac{D_h}{D_{h+1}}=(1+v+v^2)\frac{1-v^{2h+2}}{1-v^{2h+4}}.
\end{equation*}

Now let $N^{\le h}(z)$ be the generating function of Motzkin paths of height $\le h$, but where horizontal steps on level $h$ are forbidden. The system to
compute this is quite similar:
\begin{equation*}
	\begin{pmatrix}
		1-z&-z&0&\dots\\
		-z&1-z&-z&0&\dots\\
		0&-z&1-z&-z&0&\dots\\
		\ddots&\ddots&\ddots\\
		\phantom{0}&\phantom{0}&\phantom{0}&\phantom{0}&-z&\boldsymbol{1}
	\end{pmatrix}
	\begin{pmatrix}
		g_0\\g_1\\g_2\\ \vdots\\ g_h
	\end{pmatrix}
	=\begin{pmatrix}
		1\\0\\0\\ \vdots\\ 0
	\end{pmatrix}
\end{equation*}
The only difference in the matrix is the entry in the last row, written in boldface. Let $D_n^*$ be the determinant of this system matrix, with $n$ rows and columns. Again,
\begin{equation*}
	N^{\le h}(z)=g_0=\frac{D_h^*}{D_{h+1}^*}.
\end{equation*}
Expanding the matrix along the last row, we find
\begin{equation*}
D_n^*=D_{n-1}-z^2D_{n-2}=\frac1{1-v}\frac{1-v^{2n+1}}{(1+v+v^2)^n}
\end{equation*}
and
\begin{equation*}
	N^{\le h}(z)=g_0=\frac{D_h^*}{D_{h+1}^*}=(1+v+v^2)\frac{1-v^{2h+1}}{1-v^{2h+3}}.
\end{equation*}
Now let us consider
\begin{equation*}
M^{\le h}(z)-N^{\le h}(z).
\end{equation*}
There is obviously a lot of cancellation going on. The objects which are still counted, have height $=h$, and \emph{have} horizontal steps on level $h$. That is one of the two quantities that we wanted to compute, and we get
\begin{align*}
\text{Horiz}_h(z)&=(1+v+v^2)\frac{1-v^{2h+2}}{1-v^{2h+4}}-(1+v+v^2)\frac{1-v^{2h+1}}{1-v^{2h+3}}\\&=
(1+v+v^2)(1-v^{-2})\bigg[\frac{v^{2h+4}}{1-v^{2h+4}}-\frac{v^{2h+3}}{1-v^{2h+3}}\bigg].
\end{align*}

Similarly, considering $N^{\le h}(z)-M^{\le h-1}(z)$, we find that only objects are counted that have height $=h$, and \emph{no} horizontal steps on level $h$.
Thus
\begin{align*}
	\text{No-Horiz}_h(z)&=(1+v+v^2)\bigg[\frac{1-v^{2h+1}}{1-v^{2h+3}}-\frac{1-v^{2h}}{1-v^{2h+2}}\bigg]\\&=
	(1+v+v^2)(1-v^{-2})\bigg[\frac{v^{2h+3}}{1-v^{2h+3}}-\frac{v^{2h+2}}{1-v^{2h+2}}\bigg].
\end{align*}
As a check, we get
\begin{align*}
\text{Horiz}_h(z)+\text{No-Horiz}_h(z)&=(1+v+v^2)(1-v^{-2})\bigg[\frac{v^{2h+4}}{1-v^{2h+4}}-\frac{v^{2h+2}}{1-v^{2h+2}}\bigg]
\\*&=M^{\le h}(z)-M^{\le h-1}(z)=M^{= h}(z).
\end{align*}

\section{Explicit enumerations}

All our generating functions contain the term 
\begin{equation*}
(1+v+v^2)(1-v^{-2})\frac{v^{2h+a}}{1-v^{2h+a}}=(1+v+v^2)(1-v^{-2})\sum_{k\ge1}v^{k(2h+a)}
\end{equation*}
for various values of $a$. We show how to compute the coefficient of $z^n$ in this. It will be done using contour integration. The contour is a small circle in the $z$-plane or
$v$-plane.
\begin{align*}
[z^n]&(1+v+v^2)(1-v^{-2})\frac{v^{2h+a}}{1-v^{2h+a}}
=\frac1{2\pi i}\oint \frac{dz}{z^{n+1}}(1+v+v^2)(1-v^{-2})\sum_{k\ge1}v^{k(2h+a)}\\
&=\frac1{2\pi i}\oint \frac{dv(1-v^2)}{(1+v+v^2)^2}\frac{(1+v+v^2)^{n+1}}{v^{n+1}}(1+v+v^2)(1-v^{-2})\sum_{k\ge1}v^{k(2h+a)}\\
&=-\sum_{k\ge1}[v^{n+2-k(2h+a)}](1-v^2)^2(1+v+v^2)^{n}\\
&=-\sum_{k\ge1}\bigg[\binom{n,3}{n+2-k(2h+a)}-2\binom{n,3}{n-k(2h+a)}+\binom{n,3}{n-2-k(2h+a)}\bigg].
\end{align*}
In this computation, we used the notion of \emph{trinomial} coefficients:
\begin{equation*}
\binom{n,3}{k}=[t^k](1+t+t^2)^n.
\end{equation*}

\section{The average amplitude}

Here we compute:
\begin{align*}
\sum_{h\ge0}&(2h+1)\text{Horiz}_h(z)+\sum_{h\ge1}(2h)\text{No-Horiz}_h(z)\\
&=(1+v+v^2)(1-v^{-2})\sum_{h\ge0}(2h+1)\bigg[\frac{v^{2h+4}}{1-v^{2h+4}}-\frac{v^{2h+3}}{1-v^{2h+3}}\bigg]
\\&+(1+v+v^2)(1-v^{-2})\sum_{h\ge0}(2h)\bigg[\frac{v^{2h+3}}{1-v^{2h+3}}-\frac{v^{2h+2}}{1-v^{2h+2}}\bigg]\\
&=-(1+v+v^2)(1-v^{-2})\sum_{h\ge0}\frac{v^{2h+3}}{1-v^{2h+3}}\\
&+(1+v+v^2)(1-v^{-2})\sum_{h\ge0}(2h+1)\frac{v^{2h+4}}{1-v^{2h+4}}\\&-
(1+v+v^2)(1-v^{-2})\sum_{h\ge0}(2h+2)\frac{v^{2h+4}}{1-v^{2h+4}}\\
&=-(1+v+v^2)(1-v^{-2})\sum_{h\ge0}\frac{v^{2h+3}}{1-v^{2h+3}}-(1+v+v^2)(1-v^{-2})\sum_{h\ge1}\frac{v^{2h+2}}{1-v^{2h+2}}\\
&=-(1+v+v^2)(1-v^{-2})\sum_{h\ge0}\frac{v^{2h+1}}{1-v^{2h+1}}+(1+v+v^2)(1-v^{-2})\frac{v}{1-v}
\\&-(1+v+v^2)(1-v^{-2})\sum_{h\ge1}\frac{v^{2h}}{1-v^{2h}}+(1+v+v^2)(1-v^{-2})\frac{v^{2}}{1-v^{2}}\\
&=-(1+v+v^2)(1-v^{-2})\sum_{h\ge0}\frac{v^{2h+1}}{1-v^{2h+1}}
\\&-(1+v+v^2)(1-v^{-2})\sum_{h\ge1}\frac{v^{2h}}{1-v^{2h}}-\frac{(1+2v)(1+v+v^2)}{v}\\
\\&=-(1+v+v^2)(1-v^{-2})\sum_{h\ge1}\frac{v^{h}}{1-v^{h}}-\frac{(1+2v)(1+v+v^2)}{v}.
\end{align*}

To find asymptotics from here, we need the local expansion of this generating function around $v\sim1$, or, equivalently, $z\sim \frac13$.
See \cite{FGD} for more explanations of how this method works; it also involves singularity analysis of generating functions \cite{FO}.
While this combined approach of Mellin transform and singularity analysis has been used for more than 30 years, we would like to cite
\cite{steffel}, where many technical details have been worked out in a similar instance. For example, the expansion of the sum that we need is 
derived there. We combine all the expansions:
\begin{equation*}
\Big(6(1-v)+\cdots\Big)\Big(-\frac{\log(1-v)}{1-v}+\frac{\gamma}{1-v}+\cdots\Big)-9+6(1-v)+\cdots=-6\log(1-v)+\cdots
\end{equation*}
Translating this into the $z$-world means $1-v\sim \sqrt3\sqrt{1-3z}$, 
and then
\begin{equation*}
[z^n]\Big(-6\log(1-v)\Big)\sim [z^n]\Big(-6\log\sqrt{1-3z}\Big)=3\frac{3^n}{n}=\frac{3^{n+1}}{n}.
\end{equation*}
The total number of Motzkin paths is
\begin{equation*}
[z^n]\frac{1-z-\sqrt{1-2z-3z^2}}{2z^2}\sim [z^n]\Big(3-3\sqrt3\sqrt{1-3z}\Big)\sim 3\sqrt3\frac{3^{n}}{2\sqrt\pi n^{3/2}}.
\end{equation*}
Taking quotients, we get the average amplitude of random Motzkin paths of length $n$, which is asymptotic to
\begin{equation*}
2\sqrt{\frac{\pi n}{3}}.
\end{equation*}
This is about twice as much as the average height of random Motzkin paths, which is what we expected.

Now we consider
\begin{align*}
\sum_{h\ge0}\text{No-Horiz}_h(z)&=
(1+v+v^2)(1-v^{-2})\sum_{h\ge0}\bigg[\frac{v^{2h+3}}{1-v^{2h+3}}-\frac{v^{2h+2}}{1-v^{2h+2}}\bigg]\\
&=(1+v+v^2)(1-v^{-2})\sum_{h\ge1}\bigg[\frac{v^{2h-1}}{1-v^{2h-1}}-\frac{v^{2h}}{1-v^{2h}}\bigg]
+\frac{(1+v)(1+v+v^2)}{v}\\
&=(1+v+v^2)(1-v^{-2})\sum_{h\ge1}\frac{v^{h}}{1-v^{h}}(-1)^{h-1}
+\frac{(1+v)(1+v+v^2)}{v}\\
\end{align*}
This needs to be expanded about $v=1$:
\begin{equation*}
(1+v+v^2)(1-v^{-2})\sim-6(1-v)-3(1-v)^2+\cdots;
\end{equation*}
\begin{equation*}
\frac{(1+v)(1+v+v^2)}{v}\sim6-3(1-v)+\cdots;
\end{equation*}
for the remaining sum we need the Mellin transform \cite{FGD}. Set $v=e^{-t}$, and transform
\begin{equation*}
\mathscr{M}\sum_{h\ge1}\frac{v^{h}}{1-v^{h}}(-1)^{h-1}
=\mathscr{M}\sum_{h,k\ge1}e^{-thk}(-1)^{h-1}=
\Gamma(s)\zeta(s)^2(1-2^{1-s}).
\end{equation*}
By the Mellin  inversion formula:
\begin{equation*}
\sum_{h\ge1}\frac{v^{h}}{1-v^{h}}(-1)^{h-1}=
\frac1{2\pi i}\int_{2-i\infty}^{2+i\infty}\Gamma(s)\zeta(s)^2(1-2^{1-s})t^{-s}ds.
\end{equation*}
The line of integration will be shifted to the left, and the collected residues constitute the expansion about $t=0$:
\begin{align*}
\sum_{h\ge1}\frac{v^{h}}{1-v^{h}}(-1)^{h-1}\sim \frac{\log 2}{t}-\frac14+\frac{t}{48}\sim\frac{\log2}{1-v}-\frac{\log2}{2}-\frac14+\cdots.
\end{align*}
Combining,
\begin{equation*}
\sum_{h\ge0}\text{No-Horiz}_h(z)\sim6-6\log 2-\frac32(1-v)+\cdots.
\end{equation*}
But
\begin{equation*}
\frac{1-z-\sqrt{1-2z-3z^2}}{2z^2}=1+v+v^2\sim 3-3(1-v)+\cdots.
\end{equation*}
Comparing the coefficients of $1-v$, we find that asymptotically about half of the Motzkin paths belong to the `No-horizontal' and
about half of the Motzkin paths belong to the `horizontal' class. Again, this is intuitively clear, once one sees the generating functions of both families.

\bibliographystyle{plain}


\end{document}